\documentclass[12pt, a4paper]{article}

\skip\footins=8mm plus 1mm

\usepackage{latexsym,amsmath,amsfonts,amssymb,amsthm,enumerate,hyperref}
\usepackage{titlesec}
\titleformat*{\subsection}{\large\bfseries}

\usepackage[labelsep=period]{caption}

\newtheorem{theorem}{Theorem}[section]

\newtheorem{lemma}[theorem]{Lemma}

\newtheorem{problem}[theorem]{Problem}

\def\proof{\par\noindent{\textbf{Proof.}~}}

\def\Cay{{\rm Cay}}        
 \def\SL{{\rm SL}}  \def\PGL{{\rm PGL}}

\usepackage{xcolor}
\usepackage[normalem]{ulem}

\begin{document}
\openup 0.5\jot

\title{Proof of a conjecture of Green and Liebeck on codes in symmetric groups}

\author{\renewcommand{\thefootnote}{\arabic{footnote}}Teng Fang\footnotemark[1] , Jinbao Li\footnotemark[1]}

\footnotetext[1]{Department of Mathematics, Suqian University, Jiangsu 223800, PR China}

\renewcommand{\thefootnote}{}

\footnotetext{{\em E--mail addresses}: \texttt{tfangfm@foxmail.com} (T. Fang), \texttt{leejinbao25@hotmail.com} (J. Li)}

\date{}

\maketitle

\begin{abstract}

Let $A$ and $B$ be subsets of a finite group $G$ and $r$ a positive integer. If for every $g\in G$, there are precisely $r$ pairs $(a,b)\in A\times B$ such that $g=ab$, then $B$ is called a code in $G$ with respect to $A$ and we write $r G=A\boldsymbol{\cdot}B$. If in addition $B$ is a subgroup of $G$, then we say that $B$ is a subgroup code in $G$. In this paper we resolve a conjecture by Green and Liebeck \cite[Conjecture 2.3]{Green20} on certain subgroup codes in the symmetric group $S_n$. Let $n>2k$ and let $j$ be such that $2^j\leqslant k<2^{j+1}$. Suppose that $X$ is a conjugacy class in $S_n$ containing $x$, and $Y_k$ is the subgroup $S_k\times S_{n-k}$ of $S_n$, where the factor $S_k$ permutes the subset $\{1,\ldots,k\}$ and the factor $S_{n-k}$ permutes the subset $\{k+1,\ldots,n\}$. We prove that $r S_n=X\boldsymbol{\cdot}Y_k$ for some positive integer $r$ if and only if the cycle type of $x$ has exactly one cycle of length $2^i$ for $0\leqslant i\leqslant j$ and all other cycles have length at least $k+1$. We also propose several problems concerning the existence of certain subgroup codes in a finite group $G$ with respect to a conjugation-closed subset in $G$.

\medskip\smallskip
\smallskip
{\it Key words:}\;{\em code; Cayley graph; representation; symmetric group}

%\smallskip
%{\it AMS subject classification (2020):}

\end{abstract}

\section{Introduction}
\label{sec:intro}

Let $A$ be a nonempty subset of a finite group $G$ and $r$ a positive integer. We say that $A$ {\em divides} $rG$ if $G$ contains a subset $B$ such that for every $g\in G$ there are precisely $r$ pairs $(a,b)\in A\times B$ satisfying $g=a b$, in which case we write $rG=A\boldsymbol{\cdot}B$ (if $r=1$ then we write $G=A\boldsymbol{\cdot}B$), and the subset $B$ is called a {\em code} \cite{Green20,HXZ18,Terada04} in $G$ with respect to $A$. This definition of divisibility and code has corresponding interpretations in graph theory and coding theory, as illustrated below.

In a simple graph $\Gamma$ a {\em perfect code} \cite{Biggs73,Diaconis89,Terada04,vanLint82} is a subset $C$ of vertices such that every vertex of $\Gamma$ is at distance at most one from a unique vertex in $C$. Hence a set $C$ of vertices is a perfect code if and only if $\big\{\{a\}\cup \Gamma(a)\big\}_{a\in C}$ forms a partition of the vertex set of $\Gamma$, where $\Gamma(a)$ is the set of neighbors of $a$ in $\Gamma$. This generalizes the classical notion of a perfect $t$-error correcting code in the Cayley graph $H(n,q,t)$ \cite{Bannai77,Green20} and leads naturally to the study of perfect codes in Cayley graphs. Similarly, a subset $C$ of vertices is said to be a {\em total perfect code}~\cite{Zhou2016} in $\Gamma$ if every vertex of $\Gamma$ has exactly one neighbor in $C$. In graph theory, a perfect code is also called an \emph{efficient dominating set} \cite{DS2003} or {\em independent perfect dominating set}~\cite{Lee2001}, and a total perfect code is called an {\em efficient open dominating set}~\cite{HHS1998}. Given a finite group $G$ and an inverse-closed subset $A$ of $G$, the {\em Cayley graph} $\Cay(G, A)$ of $G$ with respect to the {\em connection set} $A$ is defined to have vertex set $G$, with an edge from $g$ to $h$ if and only if $gh^{-1} \in A$. Clearly, $\Cay(G, A)$ has no {loops} precisely when $1_G\not\in A$, where $1_G$ is the identity element of $G$. It is readily seen that a subset $B$ of $G$ is a code with respect to $A$ (i.e. $r G=A\boldsymbol{\cdot}B$ for some positive integer $r$) if and only if each vertex of $\Cay(G, A)$ has exactly $r$ neighbors in $B$. In particular, $G=A\boldsymbol{\cdot}B$ if and only if $B$ is a total perfect code in $\Cay(G, A)$.

\smallskip\smallskip
\noindent\textbf{{Remark}.} One sees that in a Cayley graph $\Cay(G, A)$, a subset $B$ of vertices is a perfect code if and only if $G=\big(A\cup\{1_G\}\big)\boldsymbol{\cdot}B$. Hence when $1_G\in A$, perfect codes and total perfect codes in $\Cay(G, A)$ coincide.

\smallskip
The study of (total) perfect codes in Cayley graphs has attracted quite a bit of attention due to its strong connections with coding theory, group theory and graph theory. See, for example, \cite{HXZ18} for a nice exposition of their relationship and the history of research development in this direction.

\subsection{More background}

\textbf{The $S_n$ problem.} Codes in the symmetric group $S_n$ with respect to some conjugation-closed subset are of particular interest as they are closely related to the representation theory of $S_n$. Actually, the study of codes in symmetric groups starts from as early as the 1960s. In \cite{Rothaus-Thompson}, Rothaus and Thompson considered the existence of perfect codes in $\Cay(S_n, T)$, where the subset $T$ consists of the identity and all transpositions in $S_n$. They proved that if $1+\frac{n(n-1)}{2}$ is divisible by a prime exceeding $\sqrt{n}+2$, then there does not exist a perfect code in $\Cay(S_n,T)$. For a general $n$, the existence of a perfect code in $\Cay(S_n,T)$ still remains an open problem. Interestingly, this problem may also be given a combinatorial and geometric setting \cite{Rothaus-Thompson}, and sometimes referred to as ``the $S_n$ problem'' \cite{Schmidt2002}. For $x$, $y\in S_n$, let $d(x, y)$ be the minimum number of transpositions needed to write $xy^{-1}$. One verifies that $(S_n, d)$ is a metric space, and that $T$ divides $S_n$ if and only if $S_n$ can be covered by disjoint closed spheres of radius one. This setting is closely related to the study of random walks on the symmetric group \cite{Diaconis89}.

\medskip
\noindent\textbf{Tilings of nonabelian groups.} Another interesting perspective for understanding codes in a finite group $G$ is the notion of tilings. If $G=A\boldsymbol{\cdot}B$, then the pair $(A,B)$ is also called a \emph{tiling} of $G$. This definition of tiling is figurative: if $(A,B)$ is a tiling of $G$, then the ``right-translations'' $Ab_1,\ldots,Ab_\ell$ of $A$ by elements in $B$ fill up the ``space'' $G$ without gaps or overlaps, so do the ``left-translations'' $a_1 B,\ldots,a_m B$ of $B$ by elements in $A$, where $A=\{a_1,\ldots,a_m\}$ and $B=\{b_1,\ldots,b_\ell\}$. Since the work of Haj\'{o}s~\cite{Hajos42} in his proof of a well-known conjecture of Minkowski, tilings of abelian groups have received extensive attention \cite{SS2009}. However, much less is known about tilings of nonabelian groups. As a natural
generalization of tilings, if $rG=A\boldsymbol{\cdot}B$ for some positive integer $r$, then $(A,B)$ is called an {\em $r$-tiling} of $G$. $r$-Tilings $(A,B)$ of $S_n$ and $\SL_2(q)$ such that $A$ is closed under conjugation are studied in~\cite{FLZ2021,Green20,Terada04}.

\subsection{Main results of the paper}

For a lack of known examples in the literature, Green and Liebeck \cite{Green20} were motivated to study subgroup codes in symmetric groups. Let $Y_k$ ($n>2k$) be the subgroup $S_k\times S_{n-k}$ of $S_n$, where the factor $S_k$ permutes the subset $\{1,\ldots,k\}$ and the factor $S_{n-k}$ permutes the subset $\{k+1,\ldots,n\}$. When is $Y_k$ a code in $S_n$ with respect to $X$, where $X = x^{S_n}$ is a conjugacy class in $S_n$ containing $x$? Green and Liebeck constructed several families of such codes, and based on some computational data they made a conjecture \cite[Conjecture 2.3]{Green20} about the existence of many further families of such codes in $S_n$ which we will resolve in this paper. We now state this conjecture as a theorem which is the main result of this paper.

\begin{theorem}
\label{main theorem}
Let $X$ and $Y_k$ be as above, and let $j$ be such that $2^j\leqslant k<2^{j+1}$. Then $r S_n=X\boldsymbol{\cdot}Y_k$ for some positive integer $r$ if and only if the cycle type of $x$ has exactly one cycle of length $2^i$ for $0\leqslant i\leqslant j$ and all other cycles have length at least $k+1$.

\end{theorem}

Noting that $Y_k$ is a Young subgroup (see Section~\ref{subsec:rep of Sn}) and also a maximal subgroup of $S_n$, it is of interest to consider whether other Young subgroups and maximal subgroups of $S_n$ can be codes with respect to some conjugation-closed subset (see Problem~\ref{prob:youngsub} and Problem~\ref{prob:maximalSn}). Moreover, for a finite group $G$ with certain properties (e.g. $G$ is simple) we propose the problem to determine which maximal subgroups of $G$ can be codes in $G$ with respect to some conjugation-closed subset (see Problem~\ref{prob:maximal}).

\subsection{Structure of the paper}

Our approach to prove Theorem \ref{main theorem} relies on representation theory of the symmetric group, and thus in section \ref{sec:prelim} we provide preliminaries of some notation, definitions and basic facts from general representation theory and representation theory of the symmetric group. Section \ref{sec:necsuf} and section \ref{sec:proof} are devoted to the proof of the main result. We shall first prove some necessary results and related lemmas in section \ref{sec:necsuf}, then we will tie together the results of previous sections to prove Theorem \ref{main theorem} in section \ref{sec:proof}. Finally, we propose several problems concerning the existence of certain subgroup codes in a finite group $G$ with respect to a conjugation-closed subset in $G$.

%$\chi^\lambda(x)$,  $\zeta^\lambda(x)$

\section{Preliminaries}
\label{sec:prelim}

In this section, we recall some notions and results from general representation theory and representation theory of $S_n$. For more details, the reader is referred to \cite{Curtis62,Fulton91,JK81,Sagann2001}.

%Readers familiar with the basics of this theory are invited to skip this section.

\subsection{General representation theory}

Let $\mathbb{C}G$ be the group algebra of a finite group $G$ over the complex field $\mathbb{C}$, and identify representations of $G$ (over $\mathbb{C}$) as $\mathbb{C}G$-modules. A submodule of the left regular representation of $G$ is exactly a left ideal of $\mathbb{C}G$, and an irreducible $\mathbb{C}$-linear representation of $G$ is equivalent to a minimal left ideal of $\mathbb{C}G$. By Maschke's theorem, each left ideal of $\mathbb{C}G$ is a direct sum of minimal left ideals of $\mathbb{C}G$. For a $\mathbb{C}G$-module $J$ and an irreducible representation $L$ of $G$, $L$ is said to be an {\em irreducible constituent} of $J$ if $J$ contains a submodule $I$ such that $L$ and $I$ are isomorphic as $\mathbb{C}G$-modules.

For a character $\chi$ of a representation of $G$, denote
\begin{equation}\label{equ:3}
c_\chi=\dfrac{\chi(\mathrm{1_G})}{|G|}\sum\limits_{g\in G}\chi(g^{-1})g\ \text{ and }\ I_\chi=(\mathbb{C}G)c_\chi,
\end{equation}
where $\chi(g)$ is the character value of $\chi$ evaluated at $g\in G$. Results in the following lemma are well known and can be read off from~\cite[\S33]{Curtis62}.

\begin{lemma}\label{lem:decomp}
Let $\chi_1,\dots,\chi_r$ be a complete set of irreducible characters of $G$. Then the following statements hold:
\begin{enumerate}[{\rm (a)}]
\item $c_{\chi_1},\dots,c_{\chi_r}$ form a $\mathbb{C}$-basis for the center of $\mathbb{C}G$.
\item For each $i\in\{1,\dots,r\}$, the element $c_{\chi_i}$ is the identity of $I_{\chi_i}$ and annihilates $I_{\chi_j}$ for all $j\in\{1,\dots,r\}\setminus\{i\}$. In particular, each $c_{\chi_i}$ is a central idempotent of $\mathbb{C}G$, and $I_{\chi_i}I_{\chi_j}=I_{\chi_j}I_{\chi_i}=0$ for distinct $i$ and $j$ in $\{1,\dots,r\}$.
\item $1=\sum_{i=1}^rc_{\chi_i}$, and hence $\mathbb{C}G=\bigoplus_{i=1}^rI_{\chi_i}$ is a decomposition of $\mathbb{C}G$ into simple two-sided ideals with the multiplication by $c_{\chi_i}$ being the projection into $I_{\chi_i}$.
\item For each $i\in\{1,\dots,r\}$, the two-sided ideal $I_{\chi_i}$ is the sum of all the minimal left ideals of $\mathbb{C}G$ affording the character $\chi_i$.
\end{enumerate}
\end{lemma}

\smallskip
For a subset $A$ of $G$, denote
\[
\overline{A}=\sum_{a\in A}a \in\mathbb{C}G.
\]
Then one sees that $r G=A\boldsymbol{\cdot}B$ if and only if $r\overline{G}=\overline{A}\,\overline{B}$. Actually, by definition $r G=A\boldsymbol{\cdot}B$ means that for every $g\in G$ there are precisely $r$ pairs $(a,b)\in A\times B$ satisfying $g=ab$, which is equivalent to the equation $r\overline{G}=\overline{A}\,\overline{B}$ in $\mathbb{C}G$. Suppose that $X= x^{G}$ is a conjugacy class in $G$ containing $x$. Then $\overline{X}\in Z(\mathbb{C}G)$ and the multiplication by $\overline{X}$ restricted on $I_\zeta$ is a scalar multiple, where $Z(\mathbb{C} G)$ is the center of $\mathbb{C} G$ and $\zeta$ is an irreducible character of $G$. We record these results in the next lemma.
\begin{lemma}\label{lem:prod}
The following statements hold:
\begin{enumerate}[{\rm (a)}]
\item $r G=A\boldsymbol{\cdot}B$ if and only if $r\overline{G}=\overline{A}\,\overline{B}$.
\item Suppose that $X= x^{G}$ is a conjugacy class in $G$ containing $x$ and $\zeta$ is an irreducible character of $G$. Then $\overline{X}w=w\overline{X}=\dfrac{|X|\zeta(x)}{\zeta(1_{G})}w$ for any $w\in I_\zeta$.
\end{enumerate}

\end{lemma}

\proof {\rm (b)} See~\cite[p.235]{Curtis62}.  \qed

\subsection{Representation theory of the symmetric group}
\label{subsec:rep of Sn}

Let $\lambda=(\lambda_1,\ldots,\lambda_\ell)$ be a {\em partition} of $n$ (denoted by $\lambda\vdash n$), that is, $\lambda_1+\cdots+\lambda_\ell=n$ and $\lambda_1\geqslant\cdots\geqslant\lambda_\ell\geqslant1$. The {\em Young diagram} of $\lambda$ is an array of $n$ boxes (also called cells) having $\ell$ left-justified rows with row $i$ containing $\lambda_i$ boxes for $1\leqslant i\leqslant\ell$. A {\em Young tableau of shape $\lambda$} (also called a {\em $\lambda$-tableau}) is an array $t$ obtained by replacing the boxes of the Young diagram of $\lambda$ with the numbers $1,2,\ldots,n$ bijectively. For a $\lambda$-tableau $t$, the subgroup $H(t)$ of $S_n$ fixing each row of $t$ is called the {\em horizontal group} (or {\em row group}) of $t$, and similarly the subgroup $V(t)$ of $S_n$ fixing each column of $t$ is called the {\em vertical group} (or {\em column group}) of $t$. $H(t)$ is also called a {\em Young subgroup} of (shape) $\lambda$. Let $\mu=(\mu_1,\ldots,\mu_m)\vdash n$. We say that $\lambda$ {\em dominates} $\mu$, written $\lambda\unrhd \mu$, if
$$
\lambda_1+\lambda_2+\cdots+\lambda_i\geqslant \mu_1+\mu_2+\cdots+\mu_i
$$
for all $i\geqslant1$. If $i>\ell$ (respectively, $i>m$), then we take $\lambda_i$ (respectively, $\mu_i$) to be $0$.

Let $\lambda=(\lambda_1,\ldots,\lambda_\ell)$ be a partition of $n$. Two $\lambda$-tableaux $s$ and $t$ are said to be {\em row-equivalent}, written $s\sim t$, if for each $i\in\{1,2,\ldots,\ell\}$, the sets of numbers in the $i$-th row of $s$ and $t$ are equal. A {\em $\lambda$-tabloid} is a row-equivalence class of $\lambda$-tableaux, and the $\mathbb{C}S_n$-module $M^\lambda$ corresponding to $\lambda$ is the $\mathbb{C}$-vector space with the $\lambda$-tabloids as basis and the natural permutation action of $S_n$ on $\lambda$-tabloids. For a $\lambda$-tableau $t$, we denote the tabloid $\{s\,|\,s\sim t\}$ by $[t]$ and call
\[
e_t:=\sum_{\sigma\in V(t)}\mathrm{sgn}(\sigma)\sigma[t]
\]
a {\em $\lambda$-polytabloid}. Let $S^\lambda$ denote the {\em Specht module} which is the submodule of $M^\lambda$ spanned by all $\lambda$-polytabloids. Results in the following lemma can be found in~\cite[\S2.4, \S2.11 and \S3.10]{Sagann2001}.

\begin{lemma}\label{lem:specht}
The following statements hold:
\begin{enumerate}[{\rm(a)}]
\item The set of Specht modules $\{S^\lambda\}_{\lambda\vdash n}$ forms a complete set of pairwise inequivalent irreducible $\mathbb{C}S_n$-modules.
%\item For each $\lambda\vdash n$, the dimension of $S^\lambda$ equals $n!/\prod_{x\in D(\lambda)}h^\lambda(x)$, where $h^\lambda(x)$ is the hook length of the box $x$ in the $i$-th row and $j$-th coloumn of $D(\lambda)$, namely, the number of boxes in the $a$-th row and $b$-th coloumn of $D(\lambda)$ such that either $a=i$ and $b\geq j$ or $a\geq i$ and $b=j$.
\item For each $\lambda\vdash n$, the $\mathbb{C}S_n$-module $M^\lambda$ decomposes as
\[
M^\lambda\,\cong\,\bigoplus_{\mu\,\unrhd\, \lambda} K_{\mu\lambda}S^\mu,
\]
where $K_{\mu\lambda}$ is a nonnegative integer, known as the Kostka number, such that $K_{\lambda\lambda}=1$ and that $K_{\mu\lambda}\geqslant 1$ if and only if $\mu\unrhd\lambda$.
\end{enumerate}

\end{lemma}

\medskip
Both $M^\lambda$ and $S^\lambda$ have their left ideal avatars in the group algebra $\mathbb{C}S_n$, which we now describe.
Let $\lambda\vdash n$, and let $t$ be a $\lambda$-tableau. For a subgroup $H$ of $S_n$, the mapping $w\overline{H}\mapsto wH$ (this mapping is the $\mathbb{C}$-linear span of the mapping $g \overline{H}\mapsto gH$ for each $g\in S_n$) is a $\mathbb{C} S_n$-module isomorphism from $(\mathbb{C}S_n) \overline{H}$ to the permutation representation of $S_n$ by left multiplication on the set $S_n/H$ of the left cosets of $H$. Hence, $(\mathbb{C}S_n) \overline{H(t)}$ is isomorphic to $M^\lambda$ by~\cite[Theorem 2.1.12]{Sagann2001}. Denote
\[
E_t=\sum_{\pi\in V(t)}\sum_{\rho\in H(t)}\mathrm{sgn}(\pi)\pi\rho\ \text{ and }\ L^\lambda=(\mathbb{C}S_n)E_t.
\]
It turns out that $L^\lambda$ is an avatar of $S^\lambda$ in $\mathbb{C}S_n$ (see for instance~\cite[Lemma 7.1.4]{JK81}). Therefore, we have the following lemma.

\begin{lemma}\label{lem:spechtavatar}
Let $\lambda\vdash n$, and let $t$ be a $\lambda$-tableau. Then there hold the following isomorphisms of $\mathbb{C}S_n$-modules:
\begin{enumerate}[{\rm(a)}]
\item The left ideal $(\mathbb{C}S_n)\overline{H(t)}$ of $\mathbb{C}S_n$ is isomorphic to $M^\lambda$.
\item The left ideal $L^\lambda$ of $\mathbb{C}S_n$ is isomorphic to $S^\lambda$.
\end{enumerate}
\end{lemma}

\section{Necessary and sufficient conditions}
\label{sec:necsuf}

In what follows, we will write $\chi^\lambda$ for the {character} afforded by the irreducible representation $S^\lambda$. With the notation in~\eqref{equ:3}, we write
\begin{equation*}
c^\lambda=c_{\chi^\lambda}\ \text{ and }\ I^\lambda=I_{\chi^\lambda}.
\end{equation*}

Let $Y_k=S_k\times S_{n-k}$ ($n>2k$) be the subgroup of $S_n$ in Theorem \ref{main theorem}. Obviously, $Y_k$ is a Young subgroup of the partition $(n-k,k)$. By Lemma~\ref{lem:specht} and Lemma~\ref{lem:spechtavatar} we can obtain the irreducible constituents of the left ideal $(\mathbb{C}S_n) \overline{Y_k}$ of $\mathbb{C}S_n$.
\begin{lemma}
\label{lem:irrecomp}
The Specht module $S^\lambda$ is an irreducible constituent of $(\mathbb{C}S_n) \overline{Y_k}$ if and only if $\lambda\unrhd (n-k,k)$, that is, if and only if $\lambda=(n-m,m)$ for some $m\in\{0,1,2,\ldots,k\}$. Moreover, $(\mathbb{C}S_n) \overline{Y_k}=\bigoplus\limits_{\theta\unrhd (n-k,\,k)}U^\theta$, where $U^\theta\subseteq I^\theta$ is a direct sum of $K_{\theta(n-k,\,k)}$ minimal left ideals of $\mathbb{C}S_n$.

\end{lemma}

Let $X$ and $Y_k$ be as in Theorem \ref{main theorem}. Based on the above discussion, we get the following criterion to check when $r S_n=X\boldsymbol{\cdot}Y_k$ holds.
\begin{lemma}
\label{lem:criterion}
$r S_n=X\boldsymbol{\cdot}Y_k$ for some positive integer $r$ if and only if $\chi^\lambda(x)=0$ whenever $\lambda\unrhd (n-k,k)$ and $\lambda\neq (n)$.

\end{lemma}

\proof First assume that $\chi^{(n-m,\,m)}(x)=0$, $1\leqslant m\leqslant k$. By Lemma~\ref{lem:prod}\,(a) it suffices to show that $\overline{X}\,\overline{Y_k}\in I^{(n)}$ as $I^{(n)}$ is the $\mathbb{C}$-vector space spanned by $\overline{S_n}$. Let $\mu\vdash n$. If $\mu\ntrianglerighteq (n-k,k)$, then by Lemma \ref{lem:decomp} we have $c^\mu\overline{Y_k}=0$ as $\overline{Y_k}\in\bigoplus\limits_{\theta\unrhd (n-k,\,k)}I^\theta$ by Lemma~\ref{lem:irrecomp}. If $\mu\unrhd (n-k,k)$ and $\mu\neq (n)$, i.e. if $\mu=(n-m,m)$ for some $m\in\{1,2,\ldots,k\}$, then $c^\mu\overline{X}=0$ by Lemma~\ref{lem:prod}\,(b) and our assumption. Therefore, for any $\mu\neq (n)$, we have $c^\mu\overline{X}\,\overline{Y_k}=0$, which implies $\overline{X}\,\overline{Y_k}\in I^{(n)}$ by Lemma \ref{lem:decomp}.

Next assume that $r S_n=X\boldsymbol{\cdot}Y_k$, which implies $\overline{X}\,\overline{Y_k}\in I^{(n)}$. If $\chi^{\lambda}(x)\neq 0$ for some $\lambda\unrhd (n-k,k)$ and $\lambda\neq (n)$, then $c^\lambda\overline{Y_k}\neq0$ by Lemma \ref{lem:decomp} and Lemma \ref{lem:irrecomp}, and $c^\lambda \overline{X}=d_\lambda c^\lambda$ for some $d_\lambda\in\mathbb{C}\setminus\{0\}$ by Lemma \ref{lem:prod}\,(b). We have
$$
0=c^\lambda (r\overline{S_n})=c^\lambda \overline{X}\,\overline{Y_k}=d_\lambda (c^\lambda\overline{Y_k})\neq 0,
$$
which is a contradiction. \qed

\bigskip
For the evaluation of the character value $\chi^{(n-m,\,m)}(x)$ we shall take advantage of the {\em Frobenius Formula} \cite[p.49, Frobenius Formula 4.10]{Fulton91}:
\begin{lemma}
\label{lem:frobenius}
Suppose that the cycle type of $x\in S_n$ has $i_1$ $1$-cycles, $i_2$ $2$-cycles, $\ldots$, and $i_n$ $n$-cycles. Then $\chi^{(n-m,\,m)}(x)$ equals the coefficient of the monomial $x_1^{n-m+1}x_2^{m}$ in the polynomial $\displaystyle (x_1-x_2)\prod\limits_{1\leqslant \ell\leqslant n}\big(x_1^\ell+x_2^\ell\big)^{i_\ell}$.

\end{lemma}

Now we have been fully prepared to go on to prove Theorem \ref{main theorem}.

\section{Proof of Theorem \ref{main theorem}}
\label{sec:proof}

Suppose that the cycle type of $x\in S_n$ has $i_1$ $1$-cycles, $i_2$ $2$-cycles, $\ldots$, and $i_n$ $n$-cycles. Denote by $F(x_1,x_2)$ the polynomial in Lemma \ref{lem:frobenius}. By Lemma \ref{lem:criterion}, we have
\begin{equation}\label{equ2}
  r S_n=X\boldsymbol{\cdot}Y_k\;\,\Leftrightarrow\;\,\chi^{(n-k,\,k)}(x)=0 \text{\;\,and\;\,} \tilde{r} S_n=X\boldsymbol{\cdot}Y_{k-1} \text{\;\,for some integer\;\,} \tilde{r}.
\end{equation}

\noindent Hence we can use induction on $k$ to prove Theorem \ref{main theorem}. The cases where $1\leqslant k\leqslant 3$ in Theorem \ref{main theorem} have been proved in \cite{Green20}. Now suppose that $k-1\geqslant1$ and the assertion in Theorem \ref{main theorem} holds for $k-1$. We need to show that this assertion still holds for $k$.

\medskip
\textsf{(a)\, $k=2^j$.} In this case $2^{j-1}\leqslant k-1<2^j$. Then by induction $\tilde{r} S_n=X\boldsymbol{\cdot}Y_{k-1}$ if and only if the cycle type of $x$ has exactly one cycle of length $2^i$ for $0\leqslant i\leqslant j-1$ and all other cycles have length at least $k=2^j$, or equivalently,
$$
F(x_1,x_2)=(x_1-x_2)(x_1+x_2)(x_1^2+x_2^2)(x_1^{2^{2}}+x_2^{2^{2}})\,\boldsymbol{\cdots}(x_1^{2^{j-1}}+x_2^{2^{j-1}})\prod_{2^j\leqslant \ell\leqslant n}\big(x_1^\ell+x_2^\ell\big)^{i_\ell}.
$$
Let us consider the coefficient $\gamma_k$ of $x_1^{n-k+1}x_2^k$ in the above $F(x_1,x_2)$. If we pick $x_2^k$ in one of the $i_k$ factors $x_1^k + x_2^k$ of $F(x_1,x_2)$, then it will contribute $1$ in $\gamma_k$; if we pick $x_1^k$ in all these $i_k$ factors, then this will contribute $-1$ in $\gamma_k$. Therefore, $\gamma_k=i_k-1$. On the other hand, $\gamma_k=\chi^{(n-k,\,k)}(x)$ by Lemma \ref{lem:frobenius}. Hence if $\tilde{r} S_n=X\boldsymbol{\cdot}Y_{k-1}$, then $\chi^{(n-k,\,k)}(x)=i_k-1$, and $\chi^{(n-k,\,k)}(x)=0$ if and only if $i_k=1$, which implies by (\ref{equ2}) that the assertion of Theorem \ref{main theorem} holds for $k$.

\medskip
\textsf{(b)\, $2^j<k<2^{j+1}$.} In this case $2^{j}\leqslant k-1<2^{j+1}$. Then by induction $\tilde{r} S_n=X\boldsymbol{\cdot}Y_{k-1}$ if and only if the cycle type of $x$ has exactly one cycle of length $2^i$ for $0\leqslant i\leqslant j$ and all other cycles have length at least $k$, or equivalently,
$$
F(x_1,x_2)=(x_1-x_2)(x_1+x_2)(x_1^2+x_2^2)(x_1^{2^{2}}+x_2^{2^{2}})\,\boldsymbol{\cdots}(x_1^{2^{j}}+x_2^{2^{j}})\prod_{k\leqslant \ell\leqslant n}\big(x_1^\ell+x_2^\ell\big)^{i_\ell}.
$$
Since $2^j<k<2^{j+1}$, $k$ has a unique $2$-adic expansion:
$$
k= 2^j a_j + 2^{j-1}a_{j-1} + \cdots + 2a_1 + a_0,
$$
where $a_i\in\{0,1\}$ for $0\leqslant i\leqslant j$ (note that $a_j=1$). Let $p$ be the smallest integer such that $a_p=1$. Consider the coefficient $\gamma_k$ of $x_1^{n-k+1}x_2^k$ in $F(x_1,x_2)$. The $i_k$ factors $x_1^k + x_2^k$ of $F(x_1,x_2)$ contribute $i_k$ in $\gamma_k$ if $x_2^k$ in $x_1^{n-k+1}x_2^k$ is picked from one of them. Suppose that we pick $x_1^k$ in all these $i_k$ factors. Then $x_2^k$ in $x_1^{n-k+1}x_2^k$ is obtained in exactly two ways:

\smallskip\smallskip
{\rm (i)} for $i\in \{p,p+1,\ldots, j\}$ such that $a_i=1$, pick $x_{2}^{2^i}$ in $x_1^{2^{i}}+x_2^{2^{i}}$, and this will contribute $1$ in $\gamma_k$;

\smallskip
{\rm (ii)} pick $-x_2$ in $x_1-x_2$, pick $x_{2}^{2^i}$ in $x_1^{2^{i}}+x_2^{2^{i}}$ for $0\leqslant i<p$, and pick $x_{2}^{2^i}$ in $x_1^{2^{i}}+x_2^{2^{i}}$ for $i\in \{p+1,\ldots, j\}$ such that $a_i=1$. This will contribute $-1$ in $\gamma_k$.

\smallskip\smallskip
\noindent Therefore, $\gamma_k=i_k+1-1=i_k$. On the other hand, $\gamma_k=\chi^{(n-k,\,k)}(x)$ by Lemma \ref{lem:frobenius}. Hence if $\tilde{r} S_n=X\boldsymbol{\cdot}Y_{k-1}$, then $\chi^{(n-k,\,k)}(x)=i_k$, and $\chi^{(n-k,\,k)}(x)=0$ if and only if $i_k=0$, which implies by (\ref{equ2}) that the assertion of Theorem \ref{main theorem} holds for $k$.

\smallskip
Combining \textsf{(a)} and \textsf{(b)}, by induction we have completed the proof of Theorem \ref{main theorem}. \qed

\medskip\smallskip
We conclude this paper with several problems on the existence of certain subgroup codes. First, since $Y_k$ is a Young subgroup of the partition $(n-k,k)$, it is natural to ask of which shape $\lambda\vdash n$ a Young subgroup can be a code in $S_n$ with respect to some conjugation-closed subset.
\begin{problem}\label{prob:youngsub}
{\rm  Determine all possible partitions $\lambda$ such that a Young subgroup of $\lambda$ is a code in $S_n$ with respect to some conjugation-closed subset $X$, and describe the cycle types of elements in $X$.}
\end{problem}

Second, noting that $Y_k$ is a maximal subgroup of $S_n$, a natural problem is to consider whether other maximal subgroups of $S_n$ can be codes with respect to some conjugation-closed subset. Generally, we can ask the same question for a finite group $G$ with certain properties (e.g. $G$ is simple).
\begin{problem}\label{prob:maximalSn}
{\em  Determine all possible maximal subgroups of $S_n$ which can be codes in $S_n$ with respect to some conjugation-closed subset $X$, and describe the cycle types of elements in $X$.}
\end{problem}

\begin{problem}\label{prob:maximal}
{\em  Let $G$ be a finite nonabelian simple group. Determine all possible maximal subgroups of $G$ which can be codes in $G$ with respect to some conjugation-closed subset.}
\end{problem}

\bigskip\medskip
\noindent\textbf{Acknowledgements.} The authors wish to thank Professors William Yong-Chuan Chen, Xin-gui Fang and Binzhou Xia for their advice and inspiration.

\end{document}